\documentclass[12pt]{amsart}
\usepackage{
amssymb
}

\newtheorem{theorem}{Theorem}
\newtheorem{proposition}{Proposition}
\newtheorem{lemma}{Lemma}
\newtheorem{definition}{Definition}
\newtheorem{corollary}{Corollary}

\theoremstyle{remark}
\newtheorem{remark}{Remark}

\newcommand{\eps}{\varepsilon}

\renewcommand{\r}[1]{(\ref{#1})}
\newcommand{\PDO}{$\Psi$DO}
\newcommand{\be}[1]{\begin{equation}\label{#1}}
\newcommand{\ee}{\end{equation}}

\renewcommand{\d}{\mathrm{d}}

\newcommand{\bo}{\partial M}

\title[Sharp stability estimate]{A sharp stability estimate in tensor tomography}

\author[P. Stefanov]{Plamen Stefanov}
\address{Department of Mathematics, Purdue University, West Lafayette, IN 47907}
\thanks{Partly supported by NSF}

\begin{document}
\maketitle

\section{Introduction}
Let $(M,g)$ be a compact Riemannian manifold with boundary. The geodesic ray transform $I$ of symmetric 2-tensor fields $f$ is given by
\be{01}
If(\gamma) = \int f_{ij}(\gamma(s))\dot\gamma^i(s)\dot\gamma^j(s)\, \d s,
\ee
where $\gamma$ runs over the set of all geodesics with endpoints on $\bo$. All \textit{potential fields}  $dv$ given by $(dv)_{ij} = \frac12 (\nabla_i v_j + \nabla_j v_i)$ with $v=0$ on $\bo$ belong to the kernel of $I$. The ray transform $I$ is called \textit{s-injective}  if this is the only obstruction to injectivity, i.e., if $If=0$ implies that $f$ is potential. S-injectivity can only hold under certain assumptions on $(M,g)$. A natural conjecture is that it holds on \textit{simple} manifolds, see the definition below. So far it is known to be true for some classes of simple manifolds only, including generic simple manifolds, see  \cite{Sh-book, Sh-UW, Dairbekov, SU-JAMS}. 

In the cases where s-injectivity is known, there is also a stability estimate that is not sharp. In \cite{Sh-book}, it is of conditional type with a loss of a derivative, see \r{Sh} below. In \cite{SU-JAMS}, the estimate is not of conditional type but there is still a loss of a derivative, see \r{11} below. On the other hand, if $f$ is a function, or an 1-tensor (an 1-form), there is a sharp estimate, see \cite{SU-Duke}. The purpose of this paper is to prove a sharp estimate for the ray transform of 2-tensors. 

The geodesic ray transform is a linearization of the boundary distance function and plays an important role in the inverse kinematic problem (known also as boundary or lens rigidity), see e.g., \cite{Sh-book, SU-Duke, SU-AJM, SU-lens} and the references there. There, one wants to recover $(M,g)$ given the distance function on $\bo\times\bo$ or the scattering relation $\sigma :(x,\xi)\mapsto (y,\eta)$ that maps a given $x\in \bo$ and a given incident direction $\xi$  to the exit point $y$ and the exit direction $\eta$ of the geodesic issued from $(x,\xi)$.

\section{Main Results}
\begin{definition}
We say that a compact Riemannian manifold $(M,g)$ with boundary is simple if
\begin{enumerate}
\item[(a)] The boundary $\bo$ is strictly convex, i.e.,  $\langle \nabla_{\xi}\nu,\xi\rangle > 0$ for each $\xi \in T_{x}(\bo)$ where $\nu$ is the unit outward normal to the boundary.
\item[(b)] The map $\exp_{x}:\exp_{x}^{-1}M\to M$ is a diffeomorphism for each $x\in M$.
\end{enumerate}
\end{definition}

Condition (b) implies that each pair of points is connected by a unique geodesic depending smoothly on the endpoints. It also implies that $M$ is diffeomorphic to a ball, so we can work in one fixed chart only. We will fix $M$ fixed, and choose different metrics on it. If $(M,g)$ is simple, then we call $g$ simple. 

It is known \cite{Sh-book}, see also \cite{SU-JAMS}, that each symmetric 2-tensor field $f\in L^2(M)$ admits an orthogonal decomposition
\[
f=f^s +dv,
\]
where $v$ is 1-form in $H^1_0(M)$ (vanishing on $\bo$), and $f$ is divergence free, i.e., $\delta f=0$, where $(\delta f)_j = \nabla^i f_{ij}$, and $\nabla$ is the covariant derivative. The 1-form $v$ solves $\delta d v= \delta f$, $v|_{\bo}=0$. The latter is an elliptic system and the Dirichlet boundary condition is a regular one for it, see \cite{Sh-book, SU-JAMS}. 
S-injectivity then is equivalent to the following: $If=0$ implies $f^s=0$.

Set
$$
\partial_\pm SM := \left\{(x,\omega)\in TM ; \; 
x\in \bo,\,|\omega|=1,\,
\pm \langle \omega,\nu\rangle>0 \right\},
$$
where $\nu(x)$ is the outer unit normal to $\bo$ (normal w.r.t.\ $g$, of course). Here and in what follows, we denote by $\langle \omega,\nu\rangle$ the inner product of the vectors $\omega$, $\nu$, and $|\omega|$ is meant w.r.t.\ $g$. Let $\gamma_{x,\omega}(t)$ be the (unit speed) geodesic through $(x,\omega)$, defined on its maximal interval contained in $[0,\infty)$. One can  then parametrize all maximal (directed) geodesics in $M$ by points on $\partial_-SM$, and with some abuse of notation we denote $If(x,\omega) = If(\gamma_{x,\omega})$.

One of the methods to study s-injectivity of $I$ is the energy estimates method that goes back to Mukhometov \cite{Mu1, Mu77,Mu2}, see also \cite{BGerver}, where, for simple manifolds,  injectivity of $I$    acting on functions $f$ is proved. S-injectivity (injectivity up to $\d\phi$, where $\phi=0$ on $\bo$) for 1-forms $f$ is established in \cite{AnikonovR} by a modification of the same method. The case of 2-tensors is harder. Using the so-called Pestov identity, Sharafutdinov \cite{Sh-book} showed that $I$ is s-injective under an explicit  a priori bound on the positive part of the curvature of $g$, that in particular  implies simplicity of $g$ but it is not equivalent to it. This generalized earlier results on negatively curved manifolds \cite{PestovSh}. The Pestov-Sharafutdinov approach implies the following stability estimate \cite{Sh-book} under the small curvature assumption:
\begin{equation}        \label{Sh}
\|f^s\|^2_{L^2(\Omega)} \le C \Big(  \left\|j_\nu f|_{\bo}\right\|_{L^2(\bo)} 
\|I f\|_{L^2(\partial_-SM)}  +\|I f\|^2_{H^1(\partial_-SM)}  \Big), \quad \forall f\in H^1(M),
\end{equation}
where $(j_\nu f)_i = f_{ij}\nu^j$. The Sobolev spaces on the r.h.s.\  are taken with respect to the induced measure $\d \sigma(x,\omega) $. In semigeodesic local coordinates $x=(x',x^n)$, the latter is given by $\d \sigma(x,\omega) =(\det g)^{1/2}\d x^1\dots \d x^{n-1}\allowbreak \,(\det g)^{1/2}\d\sigma_x(\omega) $, where $\d\sigma_x(\omega)=\frac1{|\xi|}\big|
\sum_{j=1}^{n}(-1)^{j-1}\xi^j d\xi^1\wedge\dots\wedge\widehat{d\xi^j}\wedge\dots\wedge d\xi^{n}\big|$.

This estimate is of conditional type --- while it implies s-injectivity (under the curvature condition), it says that $f^s$ is small if $If$ is small \textit{and} we have an a priori bound on $f$. Moreover, there is a loss of one derivative: the r.h.s.\ is finite under the condition that $f\in H^1$, while in the l.h.s., we have only the $L^2$ norm of $f^s$. 

In \cite{SU-Duke, SU-JAMS, SU-AJM}, the author and G. Uhlmann studied this problem from microlocal point of view. Introduce the following measure $\d\mu(x,\omega) = |\langle \nu,\omega\rangle|\d\sigma(x,\omega)$ on $\partial_-SM$. It is easy to show that $I : L^2(M) \to L^2(\partial_-SM, \, \d\mu)$ is bounded \cite{Sh-book} (the first space is a space of tensors, actually). The \textit{normal operator} $N=I^*I$ is well defined on $L^2(M)$ then. 
We extend slightly $M$ to a larger manifold with boundary $M_1$ that is still simple so that $M_1\Supset M$. We extend tensors defined in $M$ as zero to $M_1$. Then $I : L^2(M_1) \to L^2(\partial_-SM_1, \, \d\mu)$, and one can define $I^*I$ related to $M_1$ that a priori is different from $N$. On the other hand, it is easy to see that when restricted to tensors supported in $M$, it coincides with $N$. Then $N:L^2(M)\to L^2(M_1)$  is s-injective if and only if $I$ is s-injective. 

The main results in \cite{SU-Duke, SU-JAMS} concerning the tensor tomography problem are the following. The operator $I$ is s-injective for real analytic simple metrics. Moreover, the set of simple $C^k$metrics, where $k\gg1$ is fixed, for which $I$ is s-injective,
is open and dense in $C^k$. Therefore, we get s-injectivity for a generic set of simple metrics. 
Next, for any simple $g$ for which $I$ is s-injective, one has the stability estimate 
\be{11}
\|f^s\|_{L^2(M)} \le C \| Nf \|_{\tilde H^2(M_1)}, \quad \forall f\in H^1(M).
\ee
Here $\tilde H^2(M)$ is defined as follows. To the usual $H^1(M_1)$ norm, we add a term of the kind 
\[
\sum_{\alpha<n}\|\partial_{x^\alpha}\nabla f\|_{L^2(U)} + \|x^n\partial_{x^n}\nabla f\|_{L^2(U)} ,
\]
where $(x',x^n)$ are  semigeodesic local coordinates near $\bo$, and $U$ is any fixed neighborhood of $\bo$. The constant $C$ in \r{11} can be chosen locally uniform for $g\in C^k$. 

A natural conjecture is that $I$ is injective for all simple metrics. This is still an open problem. For any simple metric however, we have an estimate of the kind \r{11} plus the term $\|Kf\|_{L^2(M)}$, where $K$ is a smoothing operator. In Proposition~\ref{pr_1} below we prove a sharper estimate of this kind. It is also known that the solenoidal tensors on the kernel of $I$ form a finitely dimensional space of smooth tensors \cite{Sh-finite, SU-Duke, SU-JAMS, Chappa, ShSkU}. This also  follows directly from the analysis below, since the inversion problem is reduced to a Fredholm one. 

Estimate \r{11} is not of conditional type anymore but there is still a loss of one derivative. Indeed, $N$ is a \PDO\ of order $-1$, and the natural norm on the r.h.s.\ of \r{11} would be the $H^1(M_1)$ one. Our main result shows that this is the case, indeed.

\begin{theorem}  \label{thm_1}
Let $g\in C^k(M)$, $k\gg1$, be a simple metric on $M$, and assume that $I$ is s-injective. 

(a) Then
\be{12}
\|f^s\|_{L^2(M)}/C \le  \|Nf\|_{H^1(M_1)} \le C\|f^s\|_{L^2(M)}
\ee
with some $C>0$. 

(b) The constant  $C$ can be chosen uniformly under a small $C^3(M)$ perturbation of $g$. 
\end{theorem}

As pointed out above, the s-injectivity assumption is generically true for simple metrics, and holds in particular for metrics close enough to analytic ones \cite{SU-JAMS} or for metrics with an explicit bound on the curvature \cite{Sh-book}.

The new ingredient of the proof is the use of Korn's inequality \cite[Collorary~5.12.3]{Taylor-book1}, see \r{K}.

In \cite{SU-lens}, the author and G.~Uhlmann considered manifolds that are not simple, with possible conjugate points, and studied the question of the s-injectivity on $I$ known on a subset $\Gamma$ of geodesics. The basic assumption is that none of the geodesics in $\Gamma$ has conjugate points, and the conormal bundles of all $\gamma\in \Gamma$ cover $T^*M$. Under that assumption, results about s-injectivity for generic simple metrics, including analytic ones are obtained. A stability estimate of the kind \r{11} is also proven there, where $N$ is modified via a smooth cut-off that restricts the geodesics to $\Gamma$. Without going into detail, we will only mention that Theorem~\ref{thm_1} generalizes to that case, i.e., 
the stability estimate in \cite{SU-lens}
can be written in the form \r{12} as well. 

Theorem~\ref{thm_1} allows us to reduce the smoothness requirement in the generic result in \cite{SU-JAMS}.
\begin{corollary}  \label{cor_1}
There exists a dense open set of simple metrics in $C^3(M)$ so that the corresponding ray transform $I$ is s-injective (and \r{12} holds). 
\end{corollary}

Note that we are not claiming that the set of all $C^3(M)$ simple metrics with an s-injective $I$ is  open. Our success with proving estimate \r{12} that implies the openness depends on our ability to show that the problem can be reduced to a Fredholm one. We do this by constructing a parametrix, and this requires certain number $k$ of derivatives of $g$, at least $k=2n+1$. Once we have \r{12}, we use the singular operator theory to perturb \r{12} near any $g_0\in C^k$ with an s-injective $I$, by $C^3$ perturbations. 

\section{Proofs.}
\subsection{Proof of Theorem~\ref{thm_1}(a)} We start with recalling some facts from \cite{SU-Duke,SU-JAMS}. We show first that
\[
(N f)^{i'j'}(x) =  2\int_{S_xM}\omega^{i'}\omega^{j'}\int_0^\infty  f_{ij}(\gamma_{x,\omega}(t)) \dot \gamma_{x,\omega}^i(t) \dot \gamma_{x,\omega}^j(t)\,\d t\, \d\sigma_x(\omega),
\]
where $f$ is supported in $M$, and we work in $M_1$. Performing a change of variables, we get the integral representation
\begin{equation}           \label{a16}
(N  f)_{kl}(x) = 
\frac2{\sqrt{\det g(x)}}
 \int \frac{f^{ij}(y)}{\rho(x,y)^{n-1}}
\frac{\partial\rho}{\partial y^i} \frac{\partial\rho}{\partial y^j}
\frac{\partial\rho}{\partial x^k} \frac{\partial\rho}{\partial x^l} 
\Big|\!\det\frac{\partial^2(\rho^2/2)}{\partial x\partial y}\,\Big| \,\d y, \quad x\in M_1,
\end{equation}
where $\rho$ is the distance function. This form of $N$ show that $N$ is a \PDO\ of order $-1$ on the interior of $M_1$. Its principal symbol is, see \cite{SU-JAMS, S-Serdica},
\be{P4.3.1}
\sigma_p(N)^{ijkl} (x,\xi) = 2\pi\int_{S_xM_1} \omega^i\omega^j \omega^k\omega^l \delta(\xi\cdot\omega)\, \d\sigma_x(\omega),
\ee
where $\xi\cdot\omega = \xi_i\omega^i$.  This formula generalizes in an obvious way to tensors of any order.  It follows now easily that $N$ is elliptic on tensors satisfying $\xi^i f_{ij}=0$ (solenoidal tensors in the Fourier representation), and vanishes on tensors of the type $\frac12(\xi_iv_j+\xi_jv_i)$ (potential tensors in the Fourier representation). This fact allows us to construct a first order \PDO\ $Q$ so that for any $f\in L^2(M)$,
\be{20a}
QNf = f^s_{M_1}+Kf 
\ee
in $M_1$, where $f^s_{M_1}$ is the solenoidal projection of $f$ (extended as zero outside $M$) in $M_1$, and $K$ is a compact operator. We can assume that the kernel of $Q$ has a support close enough to the diagonal. The need to work in $M_1$ is due to the fact that we can use the (standard) \PDO\ calculus in an open set only. 
For more details, we refer to \cite{SU-Duke, SU-JAMS}. Note that this construction needs only a finitely smooth metric $g\in C^k(M)$, $k\gg1$, that we extend to $M_1$. If we want $K$ to be infinitely smoothing, then we need $g\in C^\infty(M)$. 

The next step is to construct $f^s$, given $f^s_{M_1}$. This can be done in an explicit way as follows. Note that 
\be{21}
f^s_{M_1}=Ef^s-dw\quad \mbox{in $M_1$}, 
\ee
where $E$ is the extension as zero to $M_1\setminus M$, and $w$ solves the elliptic system
\be{22}
\delta d w= \delta Ef^s, \quad w|_{\bo_1}=0. 
\ee
The distribution $\delta Ef^s$ is supported on $\bo$, and the solution $w$ exists in $H^1_0(M_1)$, see \cite{SU-JAMS} and Lemma~\ref{elliptic} below. In particular, $w|_{\bo}\in H^{1/2}(\bo)$ is well-defined. If we know $w|_{\bo}\in H^{1/2}(\bo)$, we can recover $w$ in $M$ because $\delta dw=0$ in the interior of $M$, by \r{22}. If we recover $w$ in $M$, we recover $f^s$ as well, in terms of $f^s_{M_1} $, by  \r{21}. Our goal therefore is to recover $w|_{\bo}$ first. 

We first determine $w$ in $M_1\setminus M$, up to a smoothing term, by the relation 
\be{22m}
f^s_{M_1}=-dw \quad \mbox{in $M_1\setminus M$},
\ee
see \r{21}. Since $w=0$ on $\bo_1$, we can integrate 
the identity
\be{22a}
\frac{\d}{\d t} w_i(\gamma)\dot \gamma^i = [dw(\gamma)]_{ij} \dot \gamma^i\dot \gamma^j
\ee
along geodesics in $M_1\setminus M$ connecting points on $\bo_1$ and $\bo$ to recover $w$ on $\bo$. Let $\tau_+(x,\xi)>0$ be characterized by $\gamma_{x,\xi}(t)\in \bo_1$ for $t=\tau_+(x,\xi)$. Then we get
\[
w_i(x)\xi^i = \int_0^{\tau_+(x,\xi)} [f^s_{M_1}]_{ij}(\gamma_{x,\xi}(t)) \dot\gamma_{x,\xi}^i(t)  \dot\gamma_{x,\xi}^j(t)  \, \d t, 
\]
for any $(x,\xi)$ so that $\{\gamma_{x,\xi}(t), \; 0\le t\le \tau_+(x,\xi)\}$ does not intersect $M$. 
That also implies easily the following non-sharp estimate
\be{22b}
\|w\|_{L^2(M_1\setminus M)} \le C\|dw\|_{L^2(M)} \le   C\| f^s_{M_1} \|_{L^2(M_1\setminus M)},
\ee
see also \cite{Sh-book} for the first inequality. 
We refer to \cite{SU-Duke,SU-JAMS} for more detail. This approach provides also a constructive way to reduce the problem to a Fredholm one. For the proof of the theorem however, this is not needed. The new ingredient in this work is that we apply Korn's inequality \cite[Collorary~5.12.3]{Taylor-book1},
\be{K}
\|w\|_{H^1(M_1\setminus M)} \le  
C \left( \|dw\|_{L^2(M_1\setminus M)}  +\|w\|_{L^2(M_1\setminus M)}   \right).
\ee
This inequality is a consequence of the fact that the Neumann boundary conditions for $\delta d$ are regular ones. Apply the trace theorem, \r{K}, \r{22b}, and \r{22m} to get
\[
\|w\|_{H^{1/2}(\bo)} \le C \|w\|_{H^1(M_1\setminus M)} \le C' \|f^s_{M_1} \|_{L^2(M_1\setminus M)} . \\
\]
Now, since $w$ solves the elliptic PDE $\delta d w=0$ in $M_1$, we get
\be{el_es}
\|w\|_{H^1(M)}  \le C\|f^s_{M_1} \|_{L^2(M_1\setminus M)},
\ee
see  Lemma~\ref{elliptic} below. 
This, together with \r{21} yields,
\be{23}
\|f^s\|_{L^2(M)}  \le \| f^s_{M_1}\|_{L^2(M)}  + C\|f^s_{M_1} \|_{L^2(M_1\setminus M)} \le 
C   \left( \|Nf\|_{H^1(M_1)}  +\|Kf\|_{L^2(M_1)}   \right).
\ee
It is worth noting that without the a priori s-injectivity assumption, we got the following.

\begin{proposition}  \label{pr_1}
For any $l>0$, there exists $k>0$ so that for any simple metric $g\in C^k(M)$,
\[
\|f^s\|_{L^2(M)}  \le 
C   \left( \|Nf\|_{H^1(M_1)}  +\|f\|_{H^{-s}(M)}   \right), \quad \forall f\in L^2(M).
\]
\end{proposition}

Estimate \r{23} (or Proposition~\ref{pr_1}), together with  \cite[Proposition~5.3.1]{Taylor-book1}, implies that if $I$, and therefore, $N$ is s-injective, then there is an estimate as above, with a different $C$, with the last term missing. This completes the proof of Theorem~\ref{thm_1}(a). 

\bigskip
We return to the elliptic regularity estimate \r{el_es}. If $\delta d$ is replaced by the Laplace operator, then \r{el_es} follows from \cite[Theorem~{5.1.3}]{Taylor-book1}. If we raise the Sobolev regularity everywhere in \r{ell} below by $1$, this just follows from the fact that the Dirichlet conditions are regular for $\delta b$. 
In our case, we  follow the proof of  \cite[Theorem~{5.1.3}]{Taylor-book1} to get the following.

\begin{lemma} \label{elliptic}
Let $u\in H^{-1}(M)$,   $\alpha\in H^{1/2}(\bo)$ be  1-forms. Then the boundary value problem
\be{2e}
\delta d w=u \quad \mbox{in $M$}, \qquad w|_{\bo} =\alpha
\ee
has a unique solution $w\in H^1(M)$, and the following estimate holds
\be{ell}
\|w\|_{H^1(M)} \le C\left( \|u\|_{H^{-1}(M)} + \|\alpha\|_{H^{1/2}(\bo)}    \right)
\ee
\end{lemma}

\begin{proof}
By a standard argument, first extend $\alpha\in H^{1/2}(M)$ to $\tilde \alpha\in H^1(M)$ by means of a fixed bounded extension operator, and then study $u-\tilde\alpha$ that satisfies homogeneous boundary conditions. This shows that we can assume that $\alpha=0$, then the boundary condition is equivalent to $w\in H^1_0(M)$. 

Note first that $\|w\|_{H^1(M)}$ and $\|dw\|_{L^2(M)}$ are equivalent norms on $H_0^1(M)$, by \r{K} and the Poincar\'e type of inequality for $dw$, see the first inequality in \r{22b}. The existence part of the theorem in the case $\alpha=0$ then follows as in [\cite[Proposition~{5.1.1}]{Taylor-book1}.

To prove the stability estimate, given $w\in H_0^1(M)$, integrate by parts to get
\[
\|w\|^2_{H^1}/C\le 
\|dw\|^2_{L^2} = -(\delta d w,w) \le \|\delta dw\|_{H^{-1}} \|w\|_{H^1}.
\]
That implies \r{2e} when $\alpha=0$, and completes the proof of the lemma.
\end{proof}

\begin{remark}
Lemma~\ref{elliptic} in particular justifies the solenoidal--potential decomposition of tensors $f$ with $L^2$ only regularity, see also \cite{SU-JAMS}. Then $f=f^s+dv$ with $v\in H^1_0$ solving $\delta d v=\delta f$.
\end{remark}

\subsection{Proof of Theorem~\ref{thm_1}(b) and Corollary~\ref{cor_1}}
To prove Corollary~\ref{cor_1}, we define the set $\mathcal{G}$ of simple metrics as follows: near any real analytic simple $g$, we choose a small enough neighborhood in the $C^3(M)$ topology, so that \r{12} still holds. Clearly, this set is open and dense. It remains to prove that this can be done, which is the statement of Theorem~\ref{thm_1}(b).

We will prove a bit more. Fix a simple $g_0\in C^3(M)$ (not necessarily real analytic) and assume that \r{12} holds; in particular, the corresponding ray transform $I$ is s-injective. We will show that there exists $0<\epsilon\ll1$ so that for any other $g$ with $\|g-g_0\|_{C^3(M)}<\epsilon$, \r{12} still holds with possibly a different constant $C>0$, independent on $g_0$ and $\epsilon$.

We will apply first \cite[Proposition~4]{FSU}. There, the weighted ray transform 
\[
I_wf(\gamma) = \int w(\gamma(s),\dot\gamma(s)) f(\gamma(s))\, \d s
\]
of functions is studied. The estimate in \cite[Proposition~4]{FSU} compares two such transforms with different weights and different metrics (actually, we study more general families of curves in \cite{FSU}). In our context, we apply \cite[Proposition~4]{FSU} to each $f_{ij}\dot\gamma^i \dot\gamma^j$, before summing up by treating $\dot\gamma^i \dot\gamma^j$ as a weight. All we need to show is that the generators of the geodesic flows related to $g_0$ and $g$ are $O(\epsilon)$ close in $C^2$. This follows from our assumption $\|g-g_0\|_{C^3(M)}<\epsilon$. Then we get
\be{p4}
\left\|(N_g - N_{g_0})f\right\|_{H^1(M_1)} \le C\epsilon \|f\|_{L^2(M)}. 
\ee
To perturb the l.h.s.\ of \r{12}, we need to compare the solenoidal projections $f^s_{g_0}$ and  $f^s_g$ of $f$ related to $g_0$ and $g$, respectively. Recall that $f^s = f-d(\delta d)_{D}^{-1}\delta f$, where $(\delta d)_{D}$ is the Dirichlet realization of the elliptic operator $\delta d$. By \cite[Lemma~1]{SU-JAMS}, $(\delta d)_{D}^{-1} : H^{-1}(M)\to H^1_0(M)$ depends continuously on $g\in C^1(M)$. The same is true for $d$ and $\delta$ in the corresponding spaces. Therefore, 
\be{p5}
\left\|f^s_g - f^s_{g_0}\right\|_{L^2(M)} \le C_0\epsilon\|f\|_{L^2(M)}.
\ee
It is enough to prove \r{12} for $f$ solenoidal, w.r.t.\ $g$, so, let us assume that. Then by \r{p4}, and our assumption that \r{12} holds for $g_0$,
\[
\|f^s_{g_0}\|_{L^2(M)} \le  C\|N_{g_0}f\|_{H^1(M_1)} \le C\epsilon \|f\|_{L^2(M)}+  C\|N_{g}f\|_{H^1(M_1)}.
\]
Apply \r{p5}, where $f^s_g=f$ to get
\[
(1-C_0\epsilon)\|f \|_{L^2(M)} \le C\epsilon \|f\|_{L^2(M)}+  C\|N_{g}f\|_{H^1(M_1)}.
\]
Therefore, if $\eps\ll1$, we still have \r{12}. This completes the proof of Theorem~\ref{thm_1}(b).  Now, Corollary~\ref{cor_1} follows immediately.


\begin{thebibliography}{10}

\bibitem{AnikonovR}
Y.~E. Anikonov and V.~G. Romanov.
\newblock On uniqueness of determination of a form of first degree by its
  integrals along geodesics.
\newblock {\em J. Inverse Ill-Posed Probl.}, 5(6):487--490 (1998), 1997.

\bibitem{BGerver}
I.~N. Bernstein and M.~L. Gerver.
\newblock A problem of integral geometry for a family of geodesics and an
  inverse kinematic seismics problem.
\newblock {\em Dokl. Akad. Nauk SSSR}, 243(2):302--305, 1978.

\bibitem{Chappa}
E.~Chappa.
\newblock On the characterization of the kernel of the geodesic {X}-ray
  transform.
\newblock {\em Trans. Amer. Math. Soc.}, 358(11):4793--4807 (electronic), 2006.

\bibitem{Dairbekov}
N.~S. Dairbekov.
\newblock Integral geometry problem for nontrapping manifolds.
\newblock {\em Inverse Problems}, 22(2):431--445, 2006.

\bibitem{FSU}
B.~Frigyik, P.~Stefanov, and G.~Uhlmann.
\newblock The {X}-ray transform for a generic family of curves and weights.
\newblock {\em J. Geom. Anal.}, 18(1):81--97, 2008.

\bibitem{Mu77}
R.~G. Muhometov.
\newblock The reconstruction problem of a two-dimensional {R}iemannian metric,
  and integral geometry.
\newblock {\em Dokl. Akad. Nauk SSSR}, 232(1):32--35, 1977.

\bibitem{Mu2}
R.~G. Muhometov.
\newblock On a problem of reconstructing {R}iemannian metrics.
\newblock {\em Sibirsk. Mat. Zh.}, 22(3):119--135, 237, 1981.

\bibitem{Mu1}
R.~G. Mukhometov.
\newblock On the problem of integral geometry ({R}ussian).
\newblock {\em Math. problems of geophysics, Akad. Nauk SSSR, Sibirsk., Otdel.,
  Vychisl., Tsentr, Novosibirsk}, 6(2):212--242, 1975.

\bibitem{PestovSh}
L.~N. Pestov and V.~A. Sharafutdinov.
\newblock Integral geometry of tensor fields on a manifold of negative
  curvature.
\newblock {\em Sibirsk. Mat. Zh.}, 29(3):114--130, 221, 1988.

\bibitem{ShSkU}
V.~Sharafutdinov, M.~Skokan, and G.~Uhlmann.
\newblock Regularity of ghosts in tensor tomography.
\newblock {\em J. Geom. Anal.}, 15(3):499--542, 2005.

\bibitem{Sh-book}
V.~A. Sharafutdinov.
\newblock {\em Integral geometry of tensor fields}.
\newblock Inverse and Ill-posed Problems Series. VSP, Utrecht, 1994.

\bibitem{Sh-finite}
V.~A. Sharafutdinov.
\newblock A finiteness theorem for the ray transform on a {R}iemannian
  manifold.
\newblock {\em Dokl. Akad. Nauk}, 355(2):167--169, 1997.

\bibitem{Sh-UW}
V.~A. Sharafutdinov.
\newblock Ray transform on Riemannian manifolds, lecture notes,
  {U}{W}--{S}eattle.
\newblock available at:
  http://www.ima.umn.edu/talks/workshops/7-16-27.2001/sharafutdinov/, 1999.

\bibitem{S-Serdica}
P.~Stefanov.
\newblock Microlocal approach to tensor tomography and boundary and lens
  rigidity.
\newblock {\em Serdica Math. J.}, 34(1):67--112, 2008.

\bibitem{SU-Duke}
P.~Stefanov and G.~Uhlmann.
\newblock Stability estimates for the {X}-ray transform of tensor fields and
  boundary rigidity.
\newblock {\em Duke Math. J.}, 123(3):445--467, 2004.

\bibitem{SU-JAMS}
P.~Stefanov and G.~Uhlmann.
\newblock Boundary rigidity and stability for generic simple metrics.
\newblock {\em J. Amer. Math. Soc.}, 18(4):975--1003 (electronic), 2005.

\bibitem{SU-lens}
P.~Stefanov and G.~Uhlmann.
\newblock Local lens rigidity with incomplete data for a class of non-simple
  riemannian manifolds.
\newblock {\em submitted}, 2007.

\bibitem{SU-AJM}
P.~Stefanov and G.~Uhlmann.
\newblock Integral geometry of tensor fields on a class of non-simple
  riemannian manifolds.
\newblock {\em Amer. J. Math.}, 130(1):239--268, 2008.

\bibitem{Taylor-book1}
M.~E. Taylor.
\newblock {\em Partial differential equations. {I}}, volume 115 of {\em Applied
  Mathematical Sciences}.
\newblock Springer-Verlag, New York, 1996.
\newblock Basic theory.

\end{thebibliography}

\end{document}